\magnification=1200
\def\A{{\cal A}}

\def\C{{\bf C}}
\def\M{{\cal M}}
\def\R{{\bf R}}

\def\hal{{\vrule height 10pt width 4pt depth 0pt}}

\centerline{\bf Geometric characterizations of some classes of operators}
\centerline{\bf in C*-algebras and von Neumann algebras}
\medskip

\centerline{Charles Akemann and Nik Weaver\footnote{*}{Second author
supported by NSF grant DMS-0070634}}

\medskip

\centerline{Dedicated to Richard V. Kadison on his 75th birthday}
\bigskip
\bigskip

{\narrower{\it
\noindent We present geometric characterizations of the partial
isometries, unitaries, and invertible operators in C*-algebras and von
Neumann algebra.
\bigskip}}
\bigskip

Let $\A$ be a unital C*-algebra. A result due to Roger Smith, given
as Lemma 4.5 of [BEZ], allows one in some sense to characterize the
unitaries in $\A$ in terms of norms of $2 \times 1$ matrices over $\A$.
This led David Blecher to ask us whether the unitaries can be recovered
from the Banach space structure of $\A$ alone, without recourse to the
product and the adjoint operations.
\medskip

One cannot expect to recover the entire C*-algebra structure of $\A$
from its Banach space structure. Indeed, the identity map is always
an isometry from $\A$ onto its opposite algebra, but there exist
C*-algebras $\A$ for which $\A$ and $\A^{op}$ are not $*$-isomorphic.
(See, e.g., [P].) Thus, linearly isometric C*-algebras need not be
$*$-isomorphic.
\medskip

However, Kadison proved in [K, Theorem 7] that any surjective
linear isometry between unitary C*-algebras can be written as a Jordan
$*$-isomomorphism followed by left multiplication by a unitary in the
range.  Paterson and Sinclair [PS] extended this result to
non-unital C*-algebras with the unitary coming from the
multiplier algebra of the range algebra.
Kadison notes in [K] that any surjective linear
isometry between unital C*-algebras takes unitaries to unitaries. Thus,
the unitaries of a C*-algebra are a Banach space invariant, which answers
Blecher's question but does not give a usable criterion for determining
which elements are unitary. By Kadison's result, one can describe the
set of unitaries of $\A$ as the orbit of the identity
(or any other unitary) under the group of isometries of $\A$ onto itself.
We asked, ``Can we get a direct characterization of unitaries without
starting with a unitary?''
\medskip

We found that we could give a fairly simple characterization of the
unitaries of $\A$ by working in the dual Banach space $\A^*$. Since only
Banach space structure is used, it is a ``geometric'' characterization.
This led to the question whether other standard classes of operators can
be similarly characterized. We were able to do this for invertible
operators and for partial isometries.
\medskip

Results of this type are not possible for other standard classes of
operators. For example, since left multiplication by any unitary of
$\A$ is an isometry from $\A$ onto itself, it follows that in general
the classes of self-adjoint, normal, and projection operators in $\A$
are not stable under isometries of $\A$. (These classes are not stable
under left multiplication by unitaries, in general.) Thus, there can be
no geometric characterizations of these classes in the above sense. In
contrast, the class of isometries in $\A$ is stable under multiplication
by any unitary; however, the identity map from $\A$ to $\A^{op}$, which
also preserves Banach space structure, takes isometries in $\A$ to
co-isometries in $\A^{op}$, so here too a geometric characterization
is impossible.
\medskip

We proceed to our results. We treat partial isometries first, as our
technique in this case is slightly more elementary.
\bigskip

\noindent {\bf Theorem 1.} {\it Let $\A$ be a C*-algebra and let
$x \in \A$, $\|x\| = 1$. Then $x$ is a partial isometry if and only if
$$\eqalign{&\{y \in \A: \hbox{ there exists }a > 0\hbox{ with }\|x + ay\| =
\|x - ay\| = 1\}\cr
= &\{y \in \A: \|x + by\| = {\rm max}(1, \|by\|)
\hbox{ for all }b \in \C\}.\cr}$$}
\medskip

\noindent {\it Proof.}
Let $X_1 = \{y \in \A:$ there exists $a > 0$ with $\|x + ay\| = \|x - ay\|
= 1\}$ and $X_2 = \{y \in \A: \|x + by\| = {\rm max}(1, \|by\|)$ for all
$b \in \C\}$. Suppose $y \in X_2$, $y \neq 0$. Then
taking $b = \pm \|y\|^{-1}$, the definition of $X_2$ yields $\|x + ay\| =
\|x - ay\| = 1$ for
$a = \|y\|^{-1}$. This shows that $X_2 \subset X_1$.
\medskip

For the reverse containment, let $\A$ be faithfully represented on a
Hilbert space $H$ and suppose $x$ is a partial isometry in $\A$. Let
$p = x^*x$ and $q = xx^*$ be the right and left support projections of
$x$, and let $y \in X_1$. For any unit vector $\xi \in {\rm ran}(p)$ we
have
$\|x(\xi)\| = 1$; so if $y(\xi) \neq 0$ then
$${\rm max}(\|x(\xi) + ay(\xi)\|, \|x(\xi) - ay(\xi)\|) > 1$$
for all $a > 0$ (because the unit ball of any Hilbert space is strictly
convex). This contradicts the assumption $y \in X_1$, so we must have
$y(\xi) = 0$. Thus, $yp = 0$, or equivalently, $y = y(1-p)$. Applying the
same argument to $x^*$ and $y^*$ yields $y = (1-q)y$, so we have $y =
(1-q)y(1-p)$. Since $x = qxp$ by [KR, 6.1.1], it follows that
$x^*y = 0 = y^*x$. Hence
$$\|x + by\|^2 = \|(x + by)^*(x + by)\| = \|x^*x + |b|^2y^*y\|=\|px^*xp
+ |b|^2(1-p)y^*y(1-p)\|,$$ thus
$$\|x + by\| = {\rm max}(\|x\|, \|by\|) = {\rm max}(1, \|by\|)$$
for all $b \in \C$, and we have shown $X_1 \subset X_2$. So if $x$ is
a partial isometry then $X_1 = X_2$.
\medskip

Now suppose $x$ is not a partial isometry. We shall construct $y \in
X_1 \setminus X_2$. Retain the faithful representation of $\A$ on $H$ and
let
$x = |x|u$ be the polar decomposition of $x$ [KR, 6.1.2]. (The absolute
value
$|x|$ of $x$ lies in $\A$, but partial isometry $u$ will generally lie
in the von Neumann algebra $B(H)$ of all bounded linear operators on $H$.)
Then
$0
\leq |x|
\leq 1$ and
$|x|$ is not a projection. Let $t$ be
a point in the spectrum of $|x|$ which is neither 0 nor 1 and let
$f: [0,1] \to [0,1]$ be a continuous function such that $f(0) = f(1) = 0$,
$f(t) > 0$, and $f(s) \leq s^{-1} - 1$ for all $s \in (0,1)$. Then
define $y = f(|x|)x$. Observe that $y \in \A$. Also, since $u$ is unitary,
$$\|x \pm y\| = \|(1 \pm f(|x|))|x|u\| = \|(1 \pm f(|x|))|x|\|.$$
By routine spectral theory arguments [T, pp. 17-21], using the facts that
$1 \in {\rm spec}(|x|)$ and
$f(1) = 0$, we get $0 \leq (1 \pm f(|x|))|x| \leq 1$,
so $\|x \pm y\| \leq 1$, hence $1 \in {\rm spec}((1 \pm f(|x|))|x|)$. Thus
$\|x
\pm y\| = 1$ and so $y \in X_1$.
\medskip

However, $y \not\in X_2$. To see this let $b = \|y\|^{-1} =
\|f(|x|)|x|\|^{-1} = \|f(s)s\|_\infty^{-1}$ (again using  [T, pp.
17-21], where $\|\cdot\|_{\infty}$ denotes the supremum norm over $s \in
[0,1]$); then
$$\|x + by\| = \|(1 + bf(|x|))|x|\| = \|(1 + bf(s))s||_\infty > 1$$
since $\|bf(s)s\|_\infty = 1$. This contradicts the equation
$\|x + by\| = {\rm max}(1, \|by\|)$, so that $y \not\in X_2$. So
if $x$ is not a partial isometry then $X_1 \neq X_2$.\hfill\hal
\bigskip

It might be worth noting that $x \in \A$, $\|x\| = 1$, is an extreme point
of ${\rm ball}(\A)$ (the unit ball of $\A$)
if and only if $X_1 = \{0\}$. The forward direction is
an immediate consequence of the definition of extreme points. For the reverse
direction, suppose $x$ is not an extreme point, and find a nonzero $y \in \A$
such that $x \pm y \in {\rm ball}(\A)$. Then
$$1 = \|x\| \leq {1\over 2}(\|x + y\| + \|x - y\|) \leq 1,$$
so that $\|x \pm y\| = 1$, and hence $y \in X_1$.
\medskip

Now if $\A$ is a C*-algebra and $x \in \A$, $\|x\| = 1$, write $S_x =
\{f \in \A^*: f(x) = \|f\| = 1\}$. The
notation comes from the fact that, if $x$ is the unit of $\A$,
then $S_x$ is the set of states. If $\M$ is a von Neumann algebra,
it is more natural and convenient to work in the predual $\M_*$; thus for
$x\in\M, \|x\| = 1$, we write $S^x =\{f \in M_*: f(x) = \|f\| = 1\}$.
\bigskip

\noindent {\bf Theorem 2.} {\it Let $\A$ be a C*-algebra and
let $x \in \A$, $\|x\| = 1$. Then $x$ is unitary if and only if $S_x$
spans $\A^*$.}
\medskip

\noindent {\it Proof.} Suppose $x$ is unitary (and hence $\A$ is unital)
and consider the map
$T: \A \to \A$ given by $T(y) = xy$. This map is a bijective isometry,
and hence so is the adjoint map $T^*: \A^* \to \A^*$. We have
$T^*f(y) = f(xy)$ for all $y \in \A$ and $f \in \A^*$. It follows that
$T^*(S_1) = T^*(S_x)$, where $1$ is the unit of $\A$ and $S_1$ is the
set of states on $\A$. Since the states span $\A^*$ [T, p. 120], and
$T^*$ is a linear isomorphism, it follows that $S_x$ spans $\A^*$.
\medskip

Conversely, suppose $S_x$ spans $\A^*$. This implies that $x$ must be an
extreme point of ${\rm ball}(\A)$ by the following argument. If $x$ is not
extreme then there exists a nonzero $y \in \A$ such that $x \pm y \in
{\rm ball}(\A)$, and then every $f \in S_x$ must satisfy
$|f(x \pm y)| \leq 1$, which together with $f(x) = 1$ implies $f(y) = 0$.
So $S_x$ cannot span $\A^*$, contradicting our hypothesis. Thus $x$ is
an extreme point as claimed, and therefore $\A$ is unital by
[S, Theorem 1.6.1]. By [K, Theorem 1] the extreme points of
the unit ball of a unital C*-algebra are partial isometries, so we
conclude that $x$ is a partial isometry.
\medskip

Suppose that $x$ is not unitary. Then either $x^*x$ or $xx^*$
is not the identity, and (since the cases are almost identical) we
may assume WLOG that $p = 1-x^*x \neq 0$.
We claim that
any $f \in S_x$ satisfies $f(p) = 0$, which implies that $S_x$ does not
span $\A^*$. Thus, let $f \in S_x$ and suppose $f(p) \neq 0$. Fix $|a| = 1$
such that $r = f(ap) > 0$, and for $t$ real consider the element $x + atp$.
Using the C*-norm condition, the facts that $p$ is a projection (hence
has norm 1), that  $x$ and $p$ have orthogonal right supports, and that
$\|x\| = |a| = 1$, we have
$$\|x + atp\|^2 = \|x^*x + |t|^2p\| \leq 1 + t^2.$$
By definition of $S_x$ and the choice of $a$, we also have
$$|f(x + atp)|^2 = (1 + rt)^2 = 1 + 2rt + r^2t^2.$$
Since $\|f\| = 1$, we get $|f(x + atp)| \leq \|x + atp\|$, and this
implies
$$1 + 2rt + r^2t^2  \leq  1 + t^2$$
for all real $t$, which is falsified by $t = r/(1-r^2)$ (unless $r = 1$,
when it is falsified by $t = 1$). This contradiction establishes that
$f(p) = 0$, as claimed.\hfill\hal
\bigskip

\noindent {\bf Theorem 3} {\it Let $\M$ be a von Neumann algebra and
let $x \in \M$, $\|x\| = 1$. Then $x$ is unitary if and only if $S^x$
spans $\M_*$.}
\medskip

\noindent {\it Proof.} The proof is the same as the proof of Theorem 2
except that we work in the predual of $\M$ instead of the dual of $\A$.
\hfill\hal
\bigskip

The characterization of invertible elements requires polar decomposition,
hence it is natural to work in a von Neumann algebra $\M$ at first.  We
shall use the characterization of unitaries in $\M$ given by Theorem 3
in our characterization of invertible elements of $\M$.
\bigskip

\noindent {\bf Theorem 4.} {\it Let $\M$ be a von Neumann algebra. Then
$x \in \M$ is invertible if and only if there exists a unitary $u \in \M$
and an $\epsilon > 0$ such that $f(x)\geq\epsilon$ for all $f \in S^u$.}
\medskip

\noindent {\it Proof.} Suppose $x$ is invertible in $\M$.  Let
$x = |x|u$ be its polar decomposition [KR, p. 401]. Since $x$ is
invertible and $u$ is a partial isometry, $u$ must be unitary. Thus $|x|
= xu^*$, so $|x|$ is positive and invertible. By spectral
theory there is an $\epsilon  > 0$ such that
$|x|\geq\epsilon 1$. Now for any $f \in S^u$ the map $y \mapsto f(yu)$ is
a state in $\M_*$ [S, Props. 1.5.2 and 1.7.8], so $f(x) = f(|x|u) \geq
\epsilon$. This proves the forward direction.
\medskip

Conversely, suppose there is a unitary $u \in \M_*$ and an $\epsilon
>   0$ such that $f(x) \geq \epsilon$ for all $f \in S^u$. For any
state $g \in \M_*$ the map $y \mapsto g(yu^*)$ belongs to
$S^u$ [S, Props. 1.5.2 and 1.7.8],
so by hypothesis $g(xu^*) \geq \epsilon$. Thus $xu^* \geq \epsilon$ by
spectral theory, which implies that $xu^*$ is invertible, hence $x$ is
invertible. \hfill\hal

\bigskip
The characterization of invertible elements in a C*-algebra $\A$ requires
that we have handy a von Neumann algebra so that we can take polar
decompositions.  This is intrinsically arranged by embedding $\A$ in
its double dual
$\A^{**}$, which has structure as a von Neumann algebra, and
the canonical embedding of $\A$ into $\A^{**}$ is a $*$-isomorphism onto
a weak* dense C*-subalgebra  [T, p. 122]. Since the predual of
$\A^{**}$ is, of course,
$\A^*$, the characterization of unitaries in
$\A^{**}$ given by Theorem 3 is intrinsic to $\A$.  The next corollary
is stated in somewhat greater generality.
\bigskip

\noindent {\bf Corollary 5.} {\it Let $\A$ be a weak* dense C*-subalgebra
of the von Neumann algebra $\M$, and let $x \in \A$. Then $x$ is
invertible in $\A$ if and only if $x$ is invertible in $\M$ if and only if
there exists a unitary $u\in\M$ and an $\epsilon > 0$ such that
$f(x)\geq\epsilon$ for all $f \in S^u$.}
\medskip

\noindent {\it Proof.} This follows immediately from Theorem 4
using the weak* continuity of multiplication in $\M$. [S, Prop. 1.7.8].
\hfill\hal
\bigskip

Now let's briefly consider adding a little bit more information and see
where it leads.  Suppose that we know that $\A$ is a unital C*-algebra and
we know which element is the unit of $\A$.  As shown in [L, proof of
Theorem 21], an element $x$ of $\A$ is self-adjoint if and only if
$||1+i\alpha x||\le 1+o(\alpha), \alpha \in \R, \alpha \to 0$.  If the
unit is assumed to be known, Lumer's criterion is geometric.  We can also
characterize self-adjointness via the condition $f(x) \in \R$ for all
$f \in S_1$. Moreover,
the involution in $\A$ can be recovered from the fact that $(x + iy)^* =
x - iy$ when $x$ and $y$ are self-adjoint. In the same
spirit we note the following characterizations of positive operators and
projections as a sample of the additional power obtained by identifying
the unit element.
\bigskip

\noindent {\bf Proposition 6.} {\it Let $\A$ be a C*-algebra with identified
unit element 1 and let $x \in \A$. The following are equivalent.

1. $x$ is positive.

2. $f(x) \ge 0$ for all $f \in S_1$.

3. $x$ is self-adjoint and $\|\|x\|\cdot 1 - x\| \leq \|x\|$.}
\bigskip
\bigskip

\noindent {\bf Proposition 7.} {\it Let $\A$ be a C*-algebra with identified
unit element 1 and let $x \in \A$. The following are equivalent.

1. $x$ is a projection.

2. $x$ is a partial isometry and $x \geq 0$.

3. There is a self-adjoint unitary $v$ in $\A$ such that $x =(1+v)/2.$}
\bigskip
\bigskip

\noindent [BEZ] D.\ P.\ Blecher, E.\ G.\ Effros, and V.\ Zarikian,
One-sided $M$-ideals and multipliers in operator spaces, I, manuscript.
\medskip

\noindent [K] R.\ V.\ Kadison, Isometries of operator algebras,
{\it Ann. of Math.\ \bf 54} (1951), 325-338.
\medskip

\noindent [KR] R.\ V.\ Kadison and J. \ Ringrose, {\it Fundamentals of the
theory of operator algebras Volume II}, Academic Press, NewYork, 1986.
\medskip

\noindent [L] G.\ Lumer, Semi-inner-product spaces, Trans. Amer. Math.
Soc. 100 (1961), 29-43.
\medskip

\noindent [P] N.\ C.\ Phillips, Continuous-trace C*-algebras not
isomorphic to their opposite algebras, manuscript.
\medskip

\noindent [PS] A.\ L. Paterson and A.\ M.\ Sinclair, Characterization of
isometries between C*-algebras, {\it J. Londaon Math. Soc.} (2)
{\bf 5} (1972), 755-761.
\medskip

\noindent [S] S.\ Sakai, {\it C*-algebras and W*-algebras}, Springer,
Heidelberg, 1971.
\medskip

\noindent [T] M.\ Takesaki, {\it Theory of operator algebras I},
Springer-Verlag, New York, 1979.

\end